\documentclass[12pt]{article}%
\usepackage{amsfonts}
\usepackage{sw20bams}
\usepackage{amsmath}
\usepackage{amssymb}
\usepackage{graphicx}%
\providecommand{\U}[1]{\protect\rule{.1in}{.1in}}
\begin{document}

\title{A Deceptively Simple Quadratic Recurrence}
\author{Steven Finch}
\date{November 7, 2024}
\maketitle

\begin{abstract}
Standard techniques for treating linear recurrences no longer apply for
quadratic recurrences. \ It is not hard to determine asymptotics for a
specific parametrized model over a wide domain of values (all $p\neq1/2$
here). \ The gap between theory and experimentation seems insurmountable,
however, at a single outlier ($p=1/2$).

\end{abstract}

\footnotetext{Copyright \copyright \ 2024 by Steven R. Finch. All rights
reserved.}Fix $0<p<1$ and define a quadratic recurrence \cite{F1-decept}
\[%
\begin{tabular}
[c]{llll}%
$a_{0}=0,$ &  & $a_{k}=(1-p)+p\,a_{k-1}^{2}$ & for $k\geq1$%
\end{tabular}
\ \ \ \
\]
which arises in the study of random Galton-Watson binary tree heights.
\ Clearly%
\[
\lim\limits_{k\rightarrow\infty}a_{k}=r=\left\{
\begin{array}
[c]{ccc}%
1 &  & \text{if }0<p\leq\dfrac{1}{2},\\
\dfrac{1-p}{p} &  & \text{if }\dfrac{1}{2}<p<1
\end{array}
\right.
\]
and $0<r\leq1$. \ We wish initially to prove that the convergence rate of
$\{a_{k}\}$ is exponential if and only if $p\neq1/2$. \ More precisely,
\[
0<\lim_{k\rightarrow\infty}\frac{r-a_{k}}{(2\,r\,p)^{k}}=r\,%
{\displaystyle\prod\limits_{j=0}^{\infty}}
\frac{r+a_{j}}{2r}<1.
\]
The case $p=1/2$ is more difficult. \ We examine Schoenfield's
\cite{Sch-decept, Sln-decept} analysis in deriving the asymptotic expansion%
\[
a_{k}\sim1-\frac{2}{k}+\frac{2\ln(k)+C}{k^{2}}-\frac{2\ln(k)^{2}%
+(2C-2)\ln(k)+(\frac{1}{2}C^{2}-C+1)}{k^{3}}+-\cdots
\]
and in calculating the constant $C=2(1.76799378...)=3.53598757...$. \ This
interesting case occurs in optimal stopping theory as well \cite{F2-decept}.

\section{Subcritical}

Assume that $0<p<1/2$. \ First, note that $0\leq a_{k}<1$ for all $k$ by
induction ($a_{k}\geq1-p>0$ is obvious; supposing $0\leq a_{k-1}<1$, we obtain
$a_{k}<(1-p)+p=1$). Now, writing $b_{k}=1-a_{k}$, we have $b_{0}=1$,
$0<b_{k}\leq1$ and
\begin{align*}
b_{k}  &  =p\left(  1-a_{k-1}^{2}\right)  =p(1-a_{k-1})(1+a_{k-1})\\
\  &  =p\,b_{k-1}(2-b_{k-1})\\
\  &  <2\,p\,b_{k-1}<(2p)^{2}b_{k-2}<(2p)^{3}b_{k-3}%
\end{align*}
thus $b_{k}<(2p)^{k}$ for all $k$. \ Observe that
\begin{align*}
b_{k}  &  =2p\,b_{k-1}\left(  1-\frac{b_{k-1}}{2}\right) \\
\  &  =(2p)^{2}b_{k-2}\left(  1-\frac{b_{k-2}}{2}\right)  \left(
1-\frac{b_{k-1}}{2}\right) \\
\  &  =(2p)^{3}b_{k-3}\left(  1-\frac{b_{k-3}}{2}\right)  \left(
1-\frac{b_{k-2}}{2}\right)  \left(  1-\frac{b_{k-1}}{2}\right) \\
\  &  =(2p)^{k}%
{\displaystyle\prod\limits_{j=0}^{k-1}}
\left(  1-\frac{b_{j}}{2}\right)
\end{align*}
hence
\begin{align*}
C  &  =\lim_{k\rightarrow\infty}\frac{1-a_{k}}{(2p)^{k}}=\lim_{k\rightarrow
\infty}\frac{b_{k}}{(2p)^{k}}=%
{\displaystyle\prod\limits_{j=0}^{\infty}}
\left(  1-\frac{b_{j}}{2}\right) \\
\  &  =%
{\displaystyle\prod\limits_{j=0}^{\infty}}
\left(  1-\frac{1-a_{j}}{2}\right)  =%
{\displaystyle\prod\limits_{j=0}^{\infty}}
\frac{1+a_{j}}{2}%
\end{align*}
exists and is nonzero since
\[%
{\displaystyle\sum\limits_{j=0}^{\infty}}
\frac{b_{j}}{2}<\frac{1}{2}\,%
{\displaystyle\sum\limits_{j=0}^{\infty}}
(2p)^{j}%
\]
converges. This completes the proof. The expression for $C$\ as an infinite
product turns out to be useful for high precision estimates of $C$, given $p$
(see Table 1). \ 

\section{Supercritical}

The following lemmata is needed for $i\geq1$:
\[
r^{i}p=r^{i-1}(1-p)
\]
and is true because $r\,p=((1-p)/p)p=1-p$. \ 

Assume that $1/2<p<1$. \ First, note that $0\leq a_{k}<r$ for all $k$ by
induction ($a_{k}\geq1-p>0$ is obvious; supposing $0\leq a_{k-1}<r$, we obtain%
\[
a_{k}<(1-p)+p\,r^{2}=p\,r+(1-p)r=r
\]
by lemmata, $i=1$ \&$\ 2$). Now, writing $b_{k}=r-a_{k}$, we have $b_{0}=r$,
$0<b_{k}\leq r$ and%
\[%
\begin{array}
[c]{lll}%
b_{k}=r-(1-p)-p\,a_{k-1}^{2} &  & \\
\;\;\;\;=r-r\,p-p\,a_{k-1}^{2} &  & \text{(by lemmata, }i=1\text{)}\\
\;\;\;\;=r(1-p)-p\,a_{k-1}^{2} &  & \\
\;\;\;\;=r^{2}p-p\,a_{k-1}^{2} &  & \text{(by lemmata, }i=2\text{)}\\
\;\;\;\;=p\left(  r^{2}-a_{k-1}^{2}\right)  &  & \\
\;\;\;\;=p(r-a_{k-1})(r+a_{k-1}) &  & \\
\;\;\;\;=p\,b_{k-1}\left[  2r-(r-a_{k-1})\right]  &  & \\
\;\;\;\;=p\,b_{k-1}(2r-b_{k-1}) &  & \\
\;\;\;\;<2\,r\,p\,b_{k-1}<(2\,r\,p)^{2}b_{k-2}<(2\,r\,p)^{3}b_{k-3} &  &
\end{array}
\]
thus $b_{k}<(2\,r\,p)^{k}b_{0}=r(2\,r\,p)^{k}$ for all $k$. Observe that
\begin{align*}
b_{k}  &  =2\,r\,p\,b_{k-1}\left(  1-\frac{b_{k-1}}{2r}\right) \\
\  &  =(2\,r\,p)^{2}b_{k-2}\left(  1-\frac{b_{k-2}}{2r}\right)  \left(
1-\frac{b_{k-1}}{2r}\right) \\
\  &  =(2\,r\,p)^{3}b_{k-3}\left(  1-\frac{b_{k-3}}{2r}\right)  \left(
1-\frac{b_{k-2}}{2r}\right)  \left(  1-\frac{b_{k-1}}{2r}\right) \\
\  &  =r(2\,r\,p)^{k}%
{\displaystyle\prod\limits_{j=0}^{k-1}}
\left(  1-\frac{b_{j}}{2r}\right)
\end{align*}
because $b_{0}=r$; hence
\begin{align*}
C  &  =\lim_{k\rightarrow\infty}\frac{r-a_{k}}{(2\,r\,p)^{k}}=\lim
_{k\rightarrow\infty}\frac{b_{k}}{(2\,r\,p)^{k}}=r\,%
{\displaystyle\prod\limits_{j=0}^{\infty}}
\left(  1-\frac{b_{j}}{2r}\right) \\
\  &  =r\,%
{\displaystyle\prod\limits_{j=0}^{\infty}}
\left(  1-\frac{r-a_{j}}{2r}\right)  =r\,%
{\displaystyle\prod\limits_{j=0}^{\infty}}
\frac{r+a_{j}}{2r}%
\end{align*}
exists and is nonzero since
\[%
{\displaystyle\sum\limits_{j=0}^{\infty}}
\frac{b_{j}}{2r}<\frac{1}{2}\,%
{\displaystyle\sum\limits_{j=0}^{\infty}}
(2\,r\,p)^{j}%
\]
converges. This completes the proof. The expression for $C$\ as an infinite
product turns out to be useful for high precision estimates of $C$, given $p$
(see Table 1).

Table 1. \textit{Numerical estimates of }$C$: \textit{no closed-form
expressions are known}\
\[%
\begin{tabular}
[c]{|l|l|l|l|l|}\hline
$p$ & $C$ &  & $p$ & $C$\\\hline
1/5 & $0.423894537869731...$ &  & 3/5 & $0.158431105979816...$\\\hline
1/4 & $0.392906852755779...$ &  & 2/3 & $0.161059687971223...$\\\hline
1/3 & $0.322119375942447...$ &  & 3/4 & $0.130968950918593...$\\\hline
2/5 & $0.237646658969724...$ &  & 4/5 & $0.105973634467432...$\\\hline
\end{tabular}
\
\]

\section{Critical}

Assume that $p=1/2$. \ While studying%
\[%
\begin{tabular}
[c]{llll}%
$a_{0}=0,$ &  & $a_{k}=\dfrac{1}{2}\left(  1+a_{k-1}^{2}\right)  $ & for
$k\geq1$%
\end{tabular}
\]
in the limit as $k\rightarrow\infty$, we shall defer to the standard texts
\cite{Msr-decept, GM-decept, Fgs-decept} just once: our starting point will be%
\[
1-a_{k}\sim\frac{2}{k+\ln(k)+\frac{1}{2}C}%
\]
or equivalently%
\[
a_{k}\sim1-\frac{2}{k}+\frac{2\ln(k)+C}{k^{2}}%
\]
for some constant $C$. \ On the basis of numerical experimentation,
Schoenfield \cite{Sch-decept, Sln-decept} hypothesized that the next terms of
the asymptotic series must be of the form%
\[
\frac{c_{3,2}\ln(k)^{2}+c_{3,1}\ln(k)+c_{3,0}}{k^{3}}+\frac{c_{4,3}\ln
(k)^{3}+c_{4,2}\ln(k)^{2}+c_{4,1}\ln(k)+c_{4,0}}{k^{4}}+\cdots.
\]
He actually went as far as $c_{20,0}/k^{20}$, but we shall stop at
$c_{4,0}/k^{4}$ for brevity's sake. \ The challenge is to express each
coefficient $c_{i,j}$ as a polynomial in $C$. \ 

Letting $x=c_{3,2}$, $y=c_{3,1}$, $z=c_{3,0}$ for convenience, we replace
$a_{k+1}$ by%
\[
1-\frac{2}{k+1}+\frac{2\ln(k+1)+C}{(k+1)^{2}}+\frac{x\ln(k+1)^{2}+y\ln
(k+1)+z}{(k+1)^{3}}%
\]
and expand in powers of $k$ and $\ln(k)$:%
\[
-\frac{2}{k+1}\sim-\frac{2}{k}+\frac{2}{k^{2}}-\frac{2}{k^{3}}+\frac{2}{k^{4}%
}-\frac{2}{k^{5}}+-\cdots,
\]%
\[
\frac{2\ln(k+1)}{(k+1)^{2}}\sim\left(  \frac{2}{k^{2}}-\frac{4}{k^{3}}%
+\frac{6}{k^{4}}-\frac{8}{k^{5}}+-\cdots\right)  \ln(k)+\left(  \frac{2}%
{k^{3}}-\frac{5}{k^{4}}+\frac{26}{3k^{5}}-+\cdots\right)  ,
\]%
\[
\frac{C}{(k+1)^{2}}\sim\frac{C}{k^{2}}-\frac{2C}{k^{3}}+\frac{3C}{k^{4}}%
-\frac{4C}{k^{5}}+-\cdots,
\]%
\[
\frac{x\ln(k+1)^{2}}{(k+1)^{3}}\sim x\left(  \frac{1}{k^{3}}-\frac{3}{k^{4}%
}+\frac{6}{k^{5}}-+\cdots\right)  \ln(k)^{2}+x\left(  \frac{2}{k^{4}}-\frac
{7}{k^{5}}+-\cdots\right)  \ln(k)+x\left(  \frac{1}{k^{5}}-+\cdots\right)  ,
\]%
\[
\frac{y\ln(k+1)}{(k+1)^{3}}\sim y\left(  \frac{1}{k^{3}}-\frac{3}{k^{4}}%
+\frac{6}{k^{5}}-+\cdots\right)  \ln(k)+y\left(  \frac{1}{k^{4}}-\frac
{7}{2k^{5}}+-\cdots\right)  ,
\]%
\[
\frac{z}{(k+1)^{3}}\sim z\left(  \frac{1}{k^{3}}-\frac{3}{k^{4}}+\frac
{6}{k^{5}}-+\cdots\right)  .
\]
Upon rearrangement, $a_{k+1}$ becomes%
\begin{align*}
&  1-3x\frac{\ln(k)^{2}}{k^{4}}+\left(  6+2x-3y\right)  \frac{\ln(k)}{k^{4}%
}+\left(  -3+3C+y-3z\right)  \frac{1}{k^{4}}+x\frac{\ln(k)^{2}}{k^{3}}\\
&  +(-4+y)\frac{\ln(k)}{k^{3}}+\left(  -2C+z\right)  \frac{1}{k^{3}}%
+2\frac{\ln(k)}{k^{2}}+\left(  2+C\right)  \frac{1}{k^{2}}-2\frac{1}{k}.
\end{align*}
Performing an analogous substitution in $a_{k}$, the expression $(1+a_{k}%
^{2})/2$ becomes%
\begin{align*}
&  1+(2-2x)\frac{\ln(k)^{2}}{k^{4}}+\left(  2C-2y\right)  \frac{\ln(k)}{k^{4}%
}+\left(  \frac{C^{2}}{2}-2z\right)  \frac{1}{k^{4}}+x\frac{\ln(k)^{2}}{k^{3}%
}\\
&  +(-4+y)\frac{\ln(k)}{k^{3}}+\left(  -2C+z\right)  \frac{1}{k^{3}}%
+2\frac{\ln(k)}{k^{2}}+\left(  2+C\right)  \frac{1}{k^{2}}-2\frac{1}{k}.
\end{align*}
Matching coefficients, we obtain%
\[%
\begin{array}
[c]{ccccc}%
-3x=2-2x &  & \text{hence} &  & c_{3,2}=x=-2;
\end{array}
\]%
\[%
\begin{array}
[c]{ccccc}%
6+2(-2)-3y=2C-2y &  & \text{hence} &  & c_{3,1}=y=-2C+2;
\end{array}
\]%
\[%
\begin{array}
[c]{ccccc}%
-3+3C+(-2C+2)-3z=\dfrac{C^{2}}{2}-2z &  & \text{hence} &  & c_{3,0}%
=z=-\dfrac{C^{2}}{2}+C-1
\end{array}
\]
as was to be shown.

Let $t=c_{4,3}$, $u=c_{4,2}$, $v=c_{4,1}$, $w=c_{4,0}$.\ \ We return to
replacing $a_{k+1}$ but with higher precision:%
\[
\frac{t\ln(k+1)^{3}}{(k+1)^{4}}\sim t\left(  \frac{1}{k^{4}}-\frac{4}{k^{5}%
}+-\cdots\right)  \ln(k)^{3}+t\left(  \frac{3}{k^{5}}-+\cdots\right)
\ln(k)^{2},
\]%
\[
\frac{u\ln(k+1)^{2}}{(k+1)^{4}}\sim u\left(  \frac{1}{k^{4}}-\frac{4}{k^{5}%
}+-\cdots\right)  \ln(k)^{2}+u\left(  \frac{2}{k^{5}}-+\cdots\right)  \ln(k),
\]%
\[
\frac{v\ln(k+1)}{(k+1)^{4}}\sim v\left(  \frac{1}{k^{4}}-\frac{4}{k^{5}%
}+-\cdots\right)  \ln(k)+v\left(  \frac{1}{k^{5}}-+\cdots\right)  ,
\]%
\[
\frac{w}{(k+1)^{4}}\sim w\left(  \frac{1}{k^{4}}-\frac{4}{k^{5}}%
+-\cdots\right)  .
\]
The new terms in the rearranged $a_{k+1}$ become%
\begin{align*}
&  -4t\frac{\ln(k)^{3}}{k^{5}}+\left(  6x+3t-4u\right)  \frac{\ln(k)^{2}%
}{k^{5}}+\left(  -8-7x+6y+2u-4v\right)  \frac{\ln(k)}{k^{5}}\\
&  +\left(  \frac{20}{3}-4C+x-\frac{7}{2}y+6z+v-4w\right)  \frac{1}{k^{5}}.
\end{align*}
Analogous substitutions in $(1+a_{k}^{2})/2$ and matching coefficients give
rise to equations%
\[
-4t=2x-2t,
\]%
\[
6x+3t-4u=Cx+2y-2u,
\]%
\[
-8-7x+6y+2u-4v=Cy+2z-2v,
\]%
\[
\frac{20}{3}-4C+x-\frac{7}{2}y+6z+v-4w=Cz-2w
\]
i.e.,%
\[%
\begin{array}
[c]{ccccc}%
-4t=2(-2)-2t &  & \text{hence} &  & c_{4,3}=t=2;
\end{array}
\]%
\[%
\begin{array}
[c]{ccccc}%
6(-2)+3(2)-4u=C(-2)+2(-2C+2)-2u &  & \text{hence} &  & c_{4,2}=u=3C-5;
\end{array}
\]%
\[
-8-7(-2)+6(-2C+2)+2(3C-5)-4v=C(-2C+2)+2\left(  -\dfrac{C^{2}}{2}+C-1\right)
-2v
\]
hence%
\[
c_{4,1}=v=\frac{3}{2}C^{2}-5C+5;
\]%
\[
\frac{20}{3}-4C+(-2)-\frac{7}{2}(-2C+2)+6\left(  -\dfrac{C^{2}}{2}+C-1\right)
+\left(  \frac{3}{2}C^{2}-5C+5\right)  -4w=C\left(  -\dfrac{C^{2}}%
{2}+C-1\right)  -2w
\]
hence$_{{}}$%
\[
c_{4,0}=w=\frac{1}{4}C^{3}-\frac{5}{4}C^{2}+\frac{5}{2}C-\frac{5}{3}.
\]
To find $c_{5,0}$ would require all expansions for $\ln(k)^{j}/k^{i}$ to order
$i=6>j$.

Schoenfield \cite{Sch-decept, Sln-decept} used an elaborate method and his
$c_{20,0}/k^{20}$ series for numerically calculating the constant $C$ to over
1000 digits. \ A less accurate method involves computing $a_{100000000}$
exactly via recursion, setting this equal to our $c_{4,0}/k^{4}$ series and
then solving:%

\[
C=3.535987572272308....
\]
In an unrelated thread, Schoenfield employed a similarly intricate procedure
to evaluate a quadratic threshold constant $\lambda=0.399524667096799...$ due
to Somos \cite{Scf-decept, Sle-decept}.

\section{Closing Words}

The recursion for $p=1/2$ appears elsewhere in disguised form. \ Letting%
\[
a_{k}=1-2\alpha_{k}%
\]
we have%
\begin{align*}
1-2\alpha_{k}  &  =\frac{1}{2}\left[  1+(1-2\alpha_{k-1})^{2}\right] \\
&  =\frac{1}{2}\left(  2-4\alpha_{k-1}+4\alpha_{k-1}^{2}\right)
=1-2\alpha_{k-1}+2\alpha_{k-1}^{2}%
\end{align*}
therefore
\[%
\begin{tabular}
[c]{llll}%
$\alpha_{0}=\dfrac{1}{2},$ &  & $\alpha_{k}=\alpha_{k-1}(1-\alpha_{k-1})$ &
for $k\geq1.$%
\end{tabular}
\ \ \ \
\]
Clearly $\alpha_{k}=(1-a_{k})/2$. \ As a consequence of the preceding,%
\[
\alpha_{k}\sim\frac{1}{k}-\frac{\ln(k)+c}{k^{2}}+\frac{\ln(k)^{2}%
+(2c-1)\ln(k)+(c^{2}-c+\frac{1}{2})}{k^{3}}-+\cdots
\]
where $c=C/2=1.767993786136154...$. \ The series $\sum_{k=0}^{\infty}%
\alpha_{k}$ diverges akin to the harmonic series \cite{FO-decept} and the
constant%
\[
s_{1}=\alpha_{0}+%
{\displaystyle\sum\limits_{k=1}^{\infty}}
\left(  \alpha_{k}-\frac{1}{k}\right)  =-1.60196478...
\]
provably exists. \ In contrast,%
\[
s_{2}=%
{\displaystyle\sum\limits_{k=0}^{\infty}}
\alpha_{k}^{2}=\frac{1}{2}%
\]
because, by induction,%
\[%
{\displaystyle\sum\limits_{k=0}^{n}}
\alpha_{k}^{2}=%
{\displaystyle\sum\limits_{k=0}^{n-1}}
\alpha_{k}^{2}+\alpha_{n}^{2}=\left(  \frac{1}{2}-\alpha_{n}\right)
+\alpha_{n}^{2}=\frac{1}{2}-\alpha_{n}(1-\alpha_{n})=\frac{1}{2}-\alpha_{n+1}%
\]
and $\alpha_{n+1}\rightarrow0$. \ Closed-form expressions for%
\[
s_{m}=%
{\displaystyle\sum\limits_{k=0}^{\infty}}
\alpha_{k}^{m}=\left\{
\begin{array}
[c]{ccc}%
0.159488853036112... &  & \text{if }m=3\\
0.068977706072225... &  & \text{if }m=4\\
0.032622409767106... &  & \text{if }m=5\\
0.015934111084642... &  & \text{if }m=6\\
0.007884618832013... &  & \text{if }m=7\\
0.003923447888623... &  & \text{if }m=8
\end{array}
\right.
\]
remain unknown. \ The constants $s_{m}$ appear in an alternative formula for
$c$. \ A technique called \textquotedblleft bootstrapping\textquotedblright%
\ is instructive here \cite{F3-decept, F4-decept}:%
\begin{align*}
\frac{1}{\alpha_{k}}  &  =\frac{1}{\alpha_{k-1}}\,\frac{1}{1-\alpha_{k-1}}\\
&  =\frac{1}{\alpha_{k-1}}\left(  1+\alpha_{k-1}+\alpha_{k-1}^{2}+\alpha
_{k-1}^{3}+\alpha_{k-1}^{4}+\cdots\right) \\
&  =\frac{1}{\alpha_{k-1}}+1+\alpha_{k-1}+\alpha_{k-1}^{2}+\alpha_{k-1}%
^{3}+\cdots\\
&  =\frac{1}{\alpha_{k-2}}+2+(\alpha_{k-1}+\alpha_{k-2})+(\alpha_{k-1}%
^{2}+\alpha_{k-2}^{2})+(\alpha_{k-1}^{3}+\alpha_{k-2}^{3})+\cdots\\
&  =\frac{1}{\alpha_{k-3}}+3+(\alpha_{k-1}+\alpha_{k-2}+\alpha_{k-3}%
)+(\alpha_{k-1}^{2}+\alpha_{k-2}^{2}+\alpha_{k-3}^{2})+(\alpha_{k-1}%
^{3}+\alpha_{k-2}^{3}+\alpha_{k-3}^{3})+\cdots\\
&  =\frac{1}{\alpha_{0}}+k+%
{\displaystyle\sum\limits_{j=0}^{k-1}}
\alpha_{j}+%
{\displaystyle\sum\limits_{j=0}^{k-1}}
\alpha_{j}^{2}+%
{\displaystyle\sum\limits_{j=0}^{k-1}}
\alpha_{j}^{3}+\cdots\\
&  \sim2+k+\left(  \ln(k)+\gamma+s_{1}\right)  +s_{2}+s_{3}+\cdots
\end{align*}
where $\gamma$ is the Euler-Mascheroni constant. \ The promised formula%
\[
c=2+\gamma+%
{\displaystyle\sum\limits_{m=1}^{\infty}}
s_{m}%
\]
unfortunately is not computationally helpful, owing to our limited
understanding of the sequence $s_{1}$, $s_{2}$, $s_{3}$, $s_{4}$, \ldots.

\section{Acknowledgements}

Robert Israel and Anthony Quas gave the simple proof that $C$ exists for
$0<p<1/2$, reproduced here from \cite{F1-decept}. \ The analogous formula
corresponding to $1/2<p<1$ is new, as far as I\ know. \ A significant portion
of the logistic recurrence asymptotics appear in \cite{IS-decept} -- I\ had
overlooked this fact -- and their claim that $\exp(c-1)\approx2.15768$ is
noteworthy. \ See \cite{F5-decept, F6-decept, F7-decept}, especially the
latter, for an extension of the critical series.  Philippe Flajolet is still
deeply missed by all who knew him.

\end{document}